\documentstyle[amssymb,12pt]{amsart}

\newtheorem{prop}{Proposition}[section]
\newtheorem{defin}{Definition}[section]
\newtheorem{lemma}{Lemma}[section]
\newtheorem{thm}{Theorem}[section]
\newtheorem{corollary}{Corollary}[section]

\theoremstyle{remark}
\newtheorem{remark}{Remark}

\begin{document}
\newcommand{\nc}{\newcommand} \nc{\on}{\operatorname}
\nc{\pa}{\partial}
\nc{\cA}{{\cal A}}\nc{\cB}{{\cal B}}\nc{\cC}{{\cal C}}
\nc{\cE}{{\cal E}}\nc{\cG}{{\cal G}}\nc{\cH}{{\cal H}}
\nc{\cX}{{\cal X}}\nc{\cR}{{\cal R}}\nc{\cL}{{\cal L}}
\nc{\cK}{{\cal K}}
\nc{\sh}{\on{sh}}\nc{\Id}{\on{Id}}\nc{\Diff}{\on{Diff}}
\nc{\ad}{\on{ad}}\nc{\Der}{\on{Der}}\nc{\End}{\on{End}}
\nc{\res}{\on{res}}\nc{\Yg}{\on{Yg}}\nc{\rat}{\on{rat}}
\nc{\Imm}{\on{Im}}\nc{\limm}{\on{lim}}\nc{\Ad}{\on{Ad}}
\nc{\Hol}{\on{Hol}}\nc{\Det}{\on{Det}}
\nc{\de}{\delta}\nc{\De}{\Delta}
\nc{\si}{\sigma}\nc{\ve}{\varepsilon}
\nc{\al}{\alpha}
\nc{\CC}{{\Bbb C}}\nc{\ZZ}{{\Bbb Z}}\nc{\NN}{{\Bbb N}}
\nc{\AAA}{{\Bbb A}}\nc{\cO}{{\cal O}} \nc{\cF}{{\cal F}}
\nc{\cW}{{\cal W}}
\nc{\la}{{\lambda}}\nc{\G}{{\frak g}}\nc{\A}{{\frak a}}
\nc{\HH}{{\frak h}}
\nc{\N}{{\frak n}}\nc{\B}{{\frak b}}
\nc{\La}{\Lambda}
\nc{\g}{\gamma}\nc{\eps}{\epsilon}\nc{\wt}{\widetilde}
\nc{\wh}{\widehat}
\nc{\bn}{\begin{equation}}\nc{\en}{\end{equation}}
\nc{\SL}{{\frak{sl}}}

%
%
%

\newcommand{\ldar}[1]{\begin{picture}(10,50)(-5,-25)
\put(0,25){\vector(0,-1){50}}
\put(5,0){\mbox{$#1$}} 
\end{picture}}

\newcommand{\lrar}[1]{\begin{picture}(50,10)(-25,-5)
\put(-25,0){\vector(1,0){50}}
\put(0,5){\makebox(0,0)[b]{\mbox{$#1$}}}
\end{picture}}

\newcommand{\luar}[1]{\begin{picture}(10,50)(-5,-25)
\put(0,-25){\vector(0,1){50}}
\put(5,0){\mbox{$#1$}}
\end{picture}}

\title[Quantum currents realization of 
elliptic quantum groups]{Quantum currents realization of the  
elliptic quantum groups $E_{\tau,\eta}({\frak{sl}}_{2})$}

\author{B. Enriquez}\address{Centre de Math\'ematiques, URA 169 du CNRS,
Ecole Polytechnique,
91128 Palai-seau, France}

\date{August 1997}
\maketitle

\begin{abstract} We review the construction by G. Felder and the author
of the realization of the elliptic quantum groups by quantum currents. 
\end{abstract}

\subsection*{}

The elliptic quantum groups were introduced by G. Felder in
\cite{F}. These are algebraic objects based on a solution $R(z,\la)$ of the 
dynamical Yang-Baxter equation. Here dynamical means that in addition to
the spectral parameter $z$, the $R$-matrix depends on a parameter $\la$,
which belongs to a product of elliptic curves, and that these parameters
undergo shifts in the various terms of the equation. 

The aim of this paper is to review the construction by G. Felder and the
author (\cite{EF}) of the realization of elliptic quantum groups
$E_{\tau,\eta}(\SL_{2})$ by
quantum current algebras. This construction relies on quasi-Hopf algebra
techniques. We introduce (sect. 2) a quantum loop algebra $U_{\hbar}\G(\tau)$
($\tau$ is the elliptic parameter, $\G = \SL_{2}$) that presents
analogies with
$E_{\tau,\eta}(\SL_{2})$. Namely, it has the property that the image in
representations of its classical $r$-matrix coincides with the classical
limit of $R(z,\la)$. $U_{\hbar}\G(\tau)$ is endowed with
``Drinfeld-type'' coproducts $\De$ and $\bar\De$ (see \cite{D-new}),
conjugated by a twist $F$ (sect. 2.3). Then our goal is 
to construct, in this algebra,
a solution of the DYBE yielding $R(z,\la)$ in finite-dimensional
representations. For that, we make use of a result of O. Babelon, D. Bernard
and E. Billey (BBB). This result extends to the dynamical situation the
theory of Drinfeld twists (see \cite{D-qH}): it states that a solution
of the so-called twisted Hopf cocycle equation, in some quasi-triangular
Hopf algebra, yields a solution of the
DYBE at the algebra level. To construct such a solution, we solve a 
factorization problem for the twist $F$ (sect. 4.3). This factorization 
in turn relies on some results on Hopf
algebra pairings within the quantum loop algebra.  
After we have obtained a DYBE solution in $U_{\hbar}\G(\tau)^{\otimes
2}$, we study its representations and construct from it $L$-operators,
that satisfy $RLL$ relations which are exactly the elliptic quantum
groups relations (sect. 4.6). 

To explain the strategy we are following, we first review (sect. 4.1)
the treatment of \cite{Rat} of the
rational analogue of our elliptic situation. In that case, the solution
of a factorization problem allows to construct a twist conjugating the
Drinfeld coproduct for Yangians with the usual coproduct. Another
derivation of this twist was earlier obtained by S. Khoroshkin and
V. Tolstoy (\cite{KT}).  

The interest of this construction lies in that it permits to 
embed the elliptic quantum group in an algebra ``with central
extension''. This allows to apply quantum Kac-Moody algebra
techniques to the study of this algebra. For example, one may expect
that the elliptic quantum KZB equations of \cite{FTV} are obtained in
terms of intertwiners or coinvariants from the algebras presented
here (such a study should be related to the apporach of \cite{ACF,ABB}
to the Ruijsenaars-Schneiders models). Another interesting application 
would be the study of the
center of the algebra  $U_{\hbar}\G(\tau)$ at the critical level,
as a  Poisson algebra, and its possible deformation as some kind of
 $\cW$-algebra, in the spirit of \cite{FrResh}. Another subject of
possible interest would be to study the connection with other types of 
elliptic quantum groups -- those arising
from the Belavin-Baxter solution or those studied in \cite{Fo}. 
Finally, it would be
interesting to find relations satisfied by a finite set of generators of
the quantum loop algebra $U_{\hbar}\G(\tau)$, or its subalgebra 
$U_{\hbar}\G_{\cO}$ (see 4.1) analogous to the
Drinfeld presentation for Yangians in terms of generators $I(a)$ and
$J(a)$. 

\section{Elliptic quantum groups}

\subsection{Function spaces associated with elliptic curves}

Let $\tau$ be a complex number of positive imaginary part, and $E$ be
the elliptic curve $\CC/\Gamma$, where $\Gamma = \ZZ +\tau\ZZ$. The
basic theta-function $\theta$ associated to $E$ is defined by the 
conditions this it is holomorphic on $\CC$, $\theta'(0) = 1$, 
the only zeroes of $\theta$ are the points of $\Gamma$,  
$\theta(z+1) = - \theta(z)$, and
$\theta ( z + \tau) = - e^{ - i\pi\tau} e^{ - 2 i \pi z }\theta(z)$. 
$\theta$ is then odd.

For $\la\in\CC$, define $L_{\la}$ as follows. If $\la$ does not belong to 
$\Gamma$, define $L_{\la}$ to be the set of holomorphic functions on 
$\CC - \Gamma$, 
$1$-periodic and such that $f(z + \tau) = e^{ - 2 i \pi \la} f(z)$. 
For $\la =0$, 
define $L_{\la}$ as the maximal isotropic subspace of $\CC((z))$ (for
the pairing $\langle f,g \rangle = \res_{0}(fgdz)$) containing all
holomorphic functions $f$ on $\CC - \Gamma$, 
$\Gamma$-periodic, such that $\oint_{a} f(z) dz =0$, where $a$
is the cycle $(i\eps,i\eps+1)$ (with $\eps$ small and $>0$). 
Finally, define $L_{\la} = e^{-2i\pi m z}L_{0}$ for $\la
= n + m\tau$.  

We then have 
\begin{equation} \label{L}
L_{\la} = \oplus_{j\ge 0} \CC \left( { {\theta ' }\over
\theta} \right)^{(j)} e^{ - 2i\pi m z}, \quad \on{if} \quad \la = n + m
\tau,  
\end{equation} 
\begin{equation} \label{Lla}
L_{\la} = \oplus_{i\ge 0}
\CC \left( {{ \theta ( \la + z )} \over
{ \theta(z)} }\right)^{(i)}, \quad \on{if} \quad \la\in \CC - \Gamma, 
\end{equation}
where we let $g' = \pa_{z} g$ and $g^{(i)} = \pa_{z}^{i}g$. 

We will set 
\begin{equation} \label{ei:la}
e_{i,\la}(z) = \left( {{ \theta ( \la + z )} \over
{ \theta(z)} }\right)^{(i)}. 
\end{equation}

\subsection{Elliptic $R$-matrix}

Let $(v_{1},v_{-1})$ be the standard basis of $\CC^{2}$ and $E_{ij}$ the
endomorphism of $\CC^{2}$ defined by $E_{ij}v_{\al} =
\de_{j\al}v_{i}$. 
Let $\hbar$ be a formal variable and set 
\begin{align} \label{R+}
R(z,\la) & = E_{11} \otimes E_{11} + E_{-1,-1} \otimes E_{-1,-1}
+ {{\theta(z)}\over{\theta(z+\hbar)}}
{{\theta(\la+\hbar)\theta(\la-\hbar)}\over{\theta(\la)^{2}}}
E_{1,1} \otimes E_{-1,-1} \\ & \nonumber
+ 
{{\theta(z)}\over{\theta(z+\hbar)}} E_{-1,-1} \otimes E_{11} 
+ {{\theta(z+\la)\theta(\hbar)}\over {\theta(z+\hbar)\theta(\la)}} 
E_{1,-1} \otimes E_{-1,1}
\\ & \nonumber
- {{\theta(z-\la)\theta(\hbar)}\over{\theta(z+\hbar)\theta(\la)}} 
E_{-1,1} \otimes E_{1,-1} . 
\end{align}

\subsection{The dynamical Yang-Baxter equation}

This matrix satisfies the equation 
\begin{align} \label{DYBE:mat}
R^{(12)}(z_{12},\la) & R^{(13)}(z_{13},\la +\hbar \bar h^{(2)})
R^{(23)}(z_{23},\la) 
\\ & \nonumber =
R^{(23)}(z_{23},\la +\hbar \bar h^{(1)}) R^{(13)}(z_{13},\la)
R^{(12)}(z_{12},\la + \hbar \bar h^{(3)}) ,  
\end{align}
where $z_{ij} = z_{i}-z_{j}$ and $z_{i},i=1,2,3$ are generic complex
numbers. We set $\bar h = E_{11} - E_{22}$, and $\bar h^{(k)}$ is the
image of $\bar h$ in the $k$th factor of $\End(\CC^{2})^{\otimes
3}$. For $i,j,k$ a permutation of $1,2,3$, and any $v\in
(\CC^{2})^{\otimes 3}$ such that $\bar h^{(k)}v = \mu v$, we set 
$R^{(ij)}(z,\la+\hbar h^{(k)}) v =R^{(ij)}(z,\la+\hbar \mu) v$.

Equation (\ref{DYBE:mat}) is called the dynamical Yang-Baxter equation.

\subsection{Elliptic quantum groups}

Set $\eta = \hbar/2$. 
The elliptic quantum group $E_{\tau,\eta}(\SL_{2})$ is defined as the
algebra generated by 
 $h$ and the
$a_{i}(\la),b_{i}(\la),c_{i}(\la),d_{i}(\la),i\ge 0,\la\in\CC - \Gamma$, 
subject to the relations 
$$
[h,a_{i}(\la)] = [h,d_{i}(\la)] = 0, \quad
[h,b_{i}(\la)] = -2b_{i}(\la),  \quad [h,c_{i}(\la)] = 2c_{i}(\la),  
$$
and if we set 
$$
a(z,\la) = \sum_{i\ge 0} a_{i}(\la) e_{i,+\hbar h/2}(z), 
\quad b(z,\la) = \sum_{i\ge 0} b_{i}(\la) e_{i,\la+\hbar (h-2)/2}(z), 
$$
$$
c(z,\la) = \sum_{i\ge 0} c_{i}(\la) e_{i,-\la-\hbar (h+2)/2}(z),  
\quad d(z,\la) = \sum_{i\ge 0} d_{i}(\la)e_{i,-\hbar h/2}(z),  
$$
and 
\begin{equation} \label{L-oper}
L(z,\la) = \pmatrix a(z,\la) & b(z,\la) \\ c(z,\la) & d(z,\la) 
\endpmatrix , 
\end{equation}
the relations
\begin{align} \label{RLL}  
 R^{(12)}(z_{1}-z_{2}, \la +\hbar h) &  
 L^{(1)}(z_{1} , \la)  L^{(2)} (z_{2},\la +\hbar h^{(1)})  
\\ & \nonumber =  L^{(2)}(z_{2}, \la)  L^{(1)}(z_{1},\la+\hbar h^{(2)})  
 R^{(12)} (\la, z_{1}-z_{2})   
\end{align} 
and 
\begin{equation} \label{Det=1}
\Det(z,\la) =
d(z-\hbar,\lambda)a(z,\lambda-\hbar)
 -
b(z-\hbar,\lambda)c(z,\lambda-\hbar)
{\theta(\lambda+\hbar h+\hbar)\over
\theta(\lambda+\hbar h)} = 1. 
\end{equation}
(These relations are made explicit in \cite{FV}.)

We use the convention that $x_{\la+\hbar h}(z) = \sum_{i\ge 0}
\pa_{\la}^{i}x_{\la}(z)(\hbar h)^{i}/i!$. 

\subsection{Another presentation}

The formulas defining the quantum groups of \cite{F} are based on an
$R$-matrix $\bar R$ 
slightly different from $R$. Let $\varphi$ be a solution
to the equation 
$$  
{{\varphi(\la -\hbar)}\over{\varphi(\la+\hbar)}} =  
{{\theta(\la)} \over{\theta(\la+\hbar)}};  
$$ 
then we have 
$$ 
\bar R(z,\la)  = 
\varphi(\la +\hbar h^{(2)}) R(z,\la) \varphi(\la+\hbar
h^{(1)})^{-1}. 
$$ 
The $L$-matrix of the elliptic quantum group based on $\bar R$ can be 
connected with that of $E_{\tau,\eta}(\SL_{2})$ using the transformation
$$ \bar L(z,\la)  = \varphi(\la +\hbar h)
L(z,\la) \varphi(\la+\hbar h^{(1)})^{-1}. 
$$ 

\section{Quantum currents algebra}

The quantum currents algebra used in \cite{EF} is an example of a family
of algebras which were introduced in \cite{HGQG}. These algebras are
associated to the data of an algebraic curve and a rational
differential. In \cite{ER}, it was shown that these algebras can be
endowed with quasi-Hopf structures, quantizing certain Manin pairs. 

\subsection{Classical structures}

\subsubsection{Manin pairs}

Let $\cK= \CC((z))$ 
be the completed local field of $E$ at its origin $0$, and $\cO =
\CC[[z]]$ 
the completed local ring at the same point. Endow $\cK$ with the scalar
product $\langle , \rangle_{\cK}$ defined by 
$$
\langle f,g \rangle_{\cK} = \res_{0}(fg dz). 
$$
Define on $\cK$ the derivation $\pa$ to be equal to $d/dz$. Then $\pa$ is
invariant w.r.t $\langle, \rangle_{\cK}$, and
$\cO$ is a maximal isotropic subring of $\cK$. 

Let us set $\A = \frak{sl}_{2}(\CC)$, and denote by $\langle ,
\rangle_{\A}$ an invariant scalar product on 
$\A$.  Let us set $\G = (\A \otimes
\cK) \oplus \CC D \oplus \CC K$; let us define on $\G$ the Lie algebra
stucture defined by the central extension of $\A \otimes \cK$
$$
c(x\otimes f,y\otimes g) = \langle x,y \rangle_{\A} \langle f,\pa g
\rangle_{\cK} K
$$
and by the derivation $[D,x\otimes f] = x \otimes \pa f$. 

Let $\G_{\cO}$ be the Lie subalgebra of $\G$ equal to 
$(\A \otimes \cO) \oplus \CC D$.  
Define $\langle , \rangle_{\A\otimes \cK}$ as the tensor product
of $\langle , \rangle_{\A}$ and $\langle , \rangle_{\cK}$, and $\langle ,
\rangle_{\G}$ as the scalar product on $\G$ defined by 
$\langle , \rangle_{\G}|_{\A \otimes \cK} = \langle , \rangle_{\A \otimes
\cK}$, $\langle D, \A \otimes \cK \rangle_{\G} = \langle K, \A \otimes \cK
\rangle_{\G} = 0$, and $\langle D,K\rangle_{\G}=1$. 
Then $\G_{\cO}$ is a maximal isotropic Lie subalgebra of $\G$. 

A maximal isotropic supplementary to $\G_{\cO}$ on $\G$ is defined by 
\begin{equation}
\G_{\la} = 
({\frak h} \otimes L_{0}) \oplus (\N_{+}\otimes L_{\la}) \oplus 
(\N_{-}\otimes L_{-\la})
\oplus \CC K . 
\end{equation}
Therefore, 
\begin{equation} \label{PM}
\G = \G_{\cO} \oplus
\G_{\la}
\end{equation}
define a Lie quasi-bialgebra structure 
on  ${\frak g}_{\cO}$, and 
(as in \cite{ER}), of double Lie quasi-bialgebra
on ${\frak g}$. Its classical 
$r$-matrix is given by the formula
$$
r_{\la} = D \otimes K + \sum_{i} {1\over 2}
h[z^{i}/i!] \otimes h[e_{i;0}] 
+
e [z^{i}/i!] \otimes f [e_{i; \la}] 
+
f [z^{i}/i!] \otimes e[e_{i; -\la}] , 
$$
because $(z^{i}/i!)_{i\ge 0},(e_{i;\la})_{i\ge 0}$ are 
dual bases of $\cO$ and $L_{\la}$; we denote $x\otimes f$ by $x[f]$; 
in other terms,
\begin{align} \label{class-r}
& r_{\la}(z,w) = 
{1\over 2}(h\otimes h){{\theta'}\over{\theta}}(z-w)
+
(e\otimes f){{\theta(z-w+\la)}\over{\theta(z-w)\theta(\la)}} 
+
(f\otimes e){{\theta(z-w-\la)}\over{\theta(z-w)\theta(-\la)}} 
\\ & \nonumber + D\otimes K . 
\end{align}

In what follows, we will set 
\begin{equation} \label{ei}
e^{i} = z^{i}/i!.
\end{equation}

\begin{remark} The expansion of $R(z-w,\la)$ in powers of $\hbar$ is  
$(\pi_{z}\otimes \pi_{w})(1+\hbar r_{\la}+\cdots)$, where $\pi_{z}$ and
$\pi_{w}$ are the $2$-dimensional evaluation represntations of $\G$ at
the points $z$ and $w$. This may be viewed as an indication that the
quantization of the Manin pair (\ref{PM}) yields a realization of the
elliptic quantum group relations. 
\end{remark}

\begin{remark} The role played by $r_{\la}(z,w)$ in the elliptic KZB
equations and the Hitchin systems are explained in \cite{FW} and
\cite{ER-Hitch}.  
\end{remark}

\subsubsection{Manin triples}

It is useful to consider the following twists of the Lie quasi-bialgebra
structures provided by (\ref{PM}). 

Let us set 
\begin{equation}
\G_{+} =  (\N_{+} \otimes \cK)
\oplus (\HH \otimes L_{0}) \oplus\CC K, \quad 
\G_{-} = (\N_{-} \otimes \cK)
\oplus (\HH \otimes L_{0}) \oplus\CC D , 
\end{equation}
and 
\begin{equation}
\bar\G_{+} =  (\N_{-} \otimes \cK)
\oplus (\HH \otimes L_{0}) \oplus\CC K
, \quad 
\bar\G_{-} = (\N_{+} \otimes \cK)
\oplus (\HH \otimes L_{0}) \oplus\CC D . 
\end{equation}

Then the decompositions 
\begin{equation} \label{TM}
\G = \G_{+} \oplus \G_{-}, \quad \G = \bar\G_{+} \oplus \bar\G_{-}
\end{equation}
are decompositions of $\G$ as direct sums of isotropic subalgebras. This
defines two Lie bialgebra structures on $\G$ Manin triples; both are
connected by a twist. Also, there is some twist connecting them with the
Lie quasi-bialgebra structure (\ref{PM}). The problem which we are going
to solve is to quantize these structures. 

\subsection{Quantum algebra}

We now present the algebra $U_{\hbar}\G(\tau)$ deforming the enveloping
algebra of $\G$. We will also denote this algebra by $A(\tau)$. 
Generators of $U_{\hbar}\G(\tau)$
are $D,K$ and the $x[\eps]$, $x=e,f,h$, $\eps\in \cK$; they are subject to
the relations 
$$
x[\al\eps] = \al x[\eps],\quad x[\eps + \eps'] = x[\eps] +
x[\eps'] , \quad \al\in\CC, \eps,\eps'\in \cK. 
$$
They serve to define the generating series 
$$
x(z) = \sum_{{i\in\ZZ}} x[\eps^{i}]\eps_{i}(z), \quad x = e,f,h,  
$$
$(\eps^{i})_{i\in\ZZ},(\eps_{i})_{i\in\ZZ}$ dual bases of $\cK$; recall
that 
$(e^{i})_{i\in\NN},(e_{i;0})_{i\in\NN}$ are dual bases of $\cO$ and $L_{0}$
and set 
$$
h^{+}(z) = \sum_{i\in\NN} h[e^{i}] e_{i;0}(z), \quad 
h^{-}(z) = \sum_{i\in\NN} h[e_{i;0}] e^{i}(z).  
$$
We will also use the series 
$$
K^{+}(z) = e^{({{q^{\pa}-q^{-\pa} }\over{2\pa}} h^{+})(z)}, \quad
K^{-}(z) = q^{h^{-}(z)}, 
$$
where $q = e^{\hbar}$. 
The relations presenting $U_{\hbar}\G(\tau)$ are then 
\begin{equation}  \label{K-K}
K\ \on{is}\ \on{central,} \quad 
[K^{+}(z) , K^{+}(w)] = [K^{-}(z) , K^{-}(w)] = 0, 
\end{equation}
\begin{align} \label{K+K-}
& \theta(z-w-\hbar) \theta(z-w+\hbar + \hbar K)
K^{+}(z)K^{-}(w) 
\\ & \nonumber 
 = \theta(z-w+\hbar) \theta(z-w-\hbar+\hbar K)
K^{-}(w) K^{+}(z) , 
\end{align}
\begin{align} \label{K+:e}
K^{+}(z)  e (w) K^{+}(z)^{-1} = { {\theta(z-w+\hbar)} \over
{\theta(z-w-\hbar)} } & e(w)
\end{align}
\begin{align} \label{K-:e}
K^{-}(z) e (w) K^{-}(z)^{-1} = { {\theta(w-z +\hbar K+\hbar)} \over
{\theta(w-z+\hbar K-\hbar)} 
} e(w),
\end{align}
\begin{equation} \label{K-f}
K^{+}(z)  f (w) K^{+}(z)^{-1} = { {\theta(w-z+\hbar)} \over
{\theta(w-z-\hbar)}} f(w),
K^{-}(z) f (w) K^{-}(z)^{-1} = { {\theta(z-w +\hbar)} \over
{\theta(z-w-\hbar)}}f(w),
\end{equation}
\begin{equation} \label{vertex-e}
\theta(z-w-\hbar) e(z) e(w) =  {\theta(z-w+\hbar)} e(w)e(z) , 
\end{equation}
\begin{equation} \label{vertex-f}
\theta(w-z-\hbar) f(z) f(w) 
= \theta(w-z+\hbar) f(w) f(z) , 
\end{equation}
\begin{equation} \label{e-f-ell}
[e(z) , f(w)] = {1 \over \hbar} \left( \delta(z,w) K^{+}(z) -
\delta(z,w-\hbar K) K^{-}(w)^{-1}\right). 
\end{equation}
Here $\delta$ denotes, as usual, the formal series $\sum_{i\in\ZZ}
z^{i}w^{-i-1}$. 

Similar relations were presented in \cite{DI}. 

\subsection{Coproducts}

The algebra $U_{\hbar}\G(\tau)$ is endowed with a Hopf structure given
by the coproduct $\Delta$ defined by 
\begin{equation} \label{Delta:K:ell}
\Delta(K^{+}(z)) = K^{+}(z) \otimes K^{+}(z), \quad
\Delta(K^{-}(z)) = K^{-}(z) \otimes K^{-}(z +\hbar K_{1}), 
\end{equation}
\begin{equation} \label{Delta:e:ell}
\Delta(e(z)) = e(z)\otimes K^{+}(z) + 1\otimes e(z),
\end{equation}
\begin{equation} \label{Delta:f:ell}
\Delta(f(z)) = f(z)\otimes 1 + K^{-}(z)^{-1} \otimes f(z+\hbar K_{1}),
\end{equation}
\begin{equation} \label{Delta:D:K:ell}
\Delta(D) = D \otimes 1 + 1\otimes D, \quad \Delta(K) = K \otimes 1 + 
1 \otimes K , 
\end{equation}
the counit $\varepsilon$, and the antipode $S$ defined by them; we set
$K_{1} = K\otimes 1, K_{2}=1\otimes K$. 

$U_{\hbar}\G(\tau)$ is also endowed with another Hopf structure given
by the coproduct $\bar\Delta$ defined by
\begin{equation} \label{bar:Delta:K:ell}
\bar\Delta(K^{+}(z)) = K^{+}(z) \otimes K^{+}(z), \quad
\bar\Delta(K^{-}(z)) = K^{-}(z) \otimes K^{-}(z + \hbar K_{1}), 
\end{equation}
\begin{equation} \label{bar:Delta:e:ell}
\bar\Delta(e(z)) = e(z -\hbar K_{2})\otimes K^{-}(z -\hbar K_{2})^{-1} 
+ 1\otimes e(z),
\end{equation}
\begin{equation} \label{bar:Delta:f:ell}
\bar\Delta(f(z)) = f(z)\otimes 1 + K^{+}(z) \otimes f(z),
\end{equation}
\begin{equation} \label{bar:Delta:D:K:ell}
\bar\Delta(D) = D \otimes 1 + 1\otimes D, \quad \bar\Delta(K) = K \otimes 1 + 
1 \otimes K , 
\end{equation}
the counit $\varepsilon$, and the antipode $\bar S$ defined by them. 

The Hopf structures associated with $\Delta$ and $\bar\Delta$ 
are connected by a twist 
\begin{equation} \label{F}
F = \exp\left(\hbar \sum_{i\in\ZZ} e[\eps_{i}] \otimes f[\eps^{i}] \right) , 
\end{equation}
where $(\eps^{i})_{\in\ZZ}$ is the basis of $\cK$ dual to 
$(\eps_{i})_{i\in\ZZ}$ w.r.t. $\langle , \rangle_{\cK}$; that is, we have 
$\bar\Delta = \Ad(F)\circ \Delta$. 

Then $F$ satisfies the cocycle equation 
\begin{equation} \label{F:cocycle}
(F\otimes 1)(\Delta\otimes 1)(F) = (1\otimes F)(1\otimes \Delta)(F)
\end{equation}
(see \cite{ER}, Prop. 3.1). 

\begin{prop}
$(U_{\hbar}\G(\tau),\De)$ and $(U_{\hbar}\G(\tau),\bar\De)$ are 
quantizations of the Manin triple structures defined by (\ref{TM}). 
The universal $R$-matrix of $(U_{\hbar}\G(\tau),\De)$ is 
$$
\cR_{\infty} =  q^{D \otimes K} q^{{1\over 2}\sum_{i\ge 0} h[e^{i}] 
\otimes h[e_{i;0}]} 
q^{\sum_{i\in\ZZ} e[\eps^{i}] \otimes f[\eps_{i}]}, 
$$
and for $(U_{\hbar}\G(\tau),\bar\De)$ it is
$F^{(21)}\cR_{\infty}F^{-1}$. 
\end{prop}

\section{The realization}

\subsection{Half-currents}

Fix a complex number $\la$ and set for $x = e,f,$ $K^+$, 
\begin{equation} \label{x+}
x^{+}_{\la}(z) = \sum_{i} x[e^{i}]e_{i;\la}(z), 
\end{equation}
and for $x=e,f,K^{-}$, 
\begin{equation} \label{x-}
x^{-}_{\la}(z) = \sum_{i} x[ e_{i;-\la} ] e^{i}(z);  
\end{equation}
recall that $(e^{i}),(e_{i;\la})$ are dual bases of $\cO$ and $L_{\la}$. 

The fields $e(z)$ and $f(z)$ are then split according to 
\begin{equation} \label{split}
e(z) = e^{+}_{\la}(z) + e^{-}_{\la}(z), \quad
f(z) = f^{+}_{-\la}(z) + f^{-}_{-\la}(z);   
\end{equation}
we call the expression $x^{\pm}_{\la}(z)$ ``half-currents''.
Let us introduce the generating series $k^{+}(z)$ and  $k^{-}(z)$,
defined by 
\begin{equation}
k^{+}(z) = e^{({{ q^{\pa} - 1}\over{2\pa}} h^{+})(z) }, 
\quad
k^{-}(z) = q^{({ 1 \over{ 1 + q^{-\pa}}} h^{-})(z) } ; 
\end{equation}
they satisfy the relations 
\begin{equation} \label{rel-k}
K^{+}(z) = k^{+}(z) k^{+}(z-\hbar), \quad
K^{-}(z) = k^{-}(z) k^{-}(z-\hbar). 
\end{equation}

\subsection{Realization}

Introduce the $L$-operators 
\begin{equation} \label{chaat}
L_{\la}^{+}(\zeta)
=
\pmatrix 1 & \theta(\hbar) f^{+}_{\la +\hbar h -\hbar}(\zeta)
\\ 0 & 1 \endpmatrix
\pmatrix k^{+}(z-\hbar) & 0 \\ 0 & k^{+}(z)^{-1} 
\endpmatrix
\pmatrix 1 & 0 \\ \hbar e^{+}_{-\la}(\zeta) & 1 \endpmatrix , 
\end{equation}
\begin{equation} \label{neila}
L^{-}_{\la}(\zeta) = \pmatrix 1 & 0 \\ \hbar e^{-}_{-\la}(\zeta-
K\hbar) & 1
\endpmatrix 
\pmatrix k^{-}(\zeta -\hbar) & 0 \\ 0 & k^{-}(\zeta)^{-1} \endpmatrix 
\pmatrix 1 & \theta(\hbar) f^{-}_{\la +\hbar h -\hbar}(\zeta ) 
\\ 0 & 1 \endpmatrix . 
\end{equation}
The main result of this paper is: 

\begin{thm}
Set 
\begin{align} \label{R-mat} 
R^{-}(z,\la) & =  E_{11} \otimes E_{11} + E_{-1,-1} \otimes E_{-1,-1} 
+ 
{{\theta(z)}\over {\theta(z -\hbar)}}  
E_{11} \otimes E_{-1,-1} 
\\ & \nonumber 
+ { {\theta(\la +\hbar) \theta(\la -\hbar)} 
\over{\theta(\la)^{2}}} 
{{\theta(z)}\over{\theta(z-\hbar)}}  
E_{-1,-1} \otimes E_{11} 
- { {\theta(z+\la) \theta(\hbar)}  
\over {\theta(z-\hbar) \theta(\la)}}
E_{1,-1} \otimes E_{-1,1} 
\\ & \nonumber  
+
{{\theta(z-\la)\theta(\hbar)} 
\over{\theta(z-\hbar)\theta(\la)}}
E_{-1,1} \otimes E_{1,-1} , 
\end{align} 
and $R^{+}(z,\la) = R(z,\la)$. 
The $L$-operators $L_{\la}^{\pm}(\zeta)$ satisfy  the relations 
\begin{equation} \label{Lpm:Lpm}
R^{\pm}(\zeta-\zeta',\la) L^{\pm(1)}_{\la+\hbar h^{(2)}} (\zeta)
L^{\pm(2)}_{\la} (\zeta')
=
L^{\pm(2)}_{\la+\hbar h^{(1)}} (\zeta') L^{\pm(1)}_{\la} (\zeta)
R^{\pm}(\zeta-\zeta',\la+\hbar h) 
\end{equation}
\begin{align} \label{L+:L-}
L^{-(1)}_{\la}(\zeta) & R^{-}(\zeta-\zeta',\la+\hbar h)
L^{+(2)}_{\la}(\zeta')
\\ & \nonumber =
L^{+(2)}_{\la+\hbar h^{(1)}} (\zeta') R^{-}(\zeta-\zeta'-K\hbar,\la)
L^{-(1)}_{\la+\hbar h^{(2)}}(\zeta) {{A(\zeta,\zeta'+K\hbar)}
\over{A(\zeta,\zeta')}},  
\end{align}
where 
\begin{equation} \label{azz'}
A(\zeta,\zeta') = \exp( 
\sum_{i\ge 0}\left({1\over
\pa}{{q^{\pa}-1}\over{q^{\pa}+1}} e^{i}\right)(\zeta)
e_{i;0}(\zeta')).
\end{equation} 
\end{thm}

This result can be viewed as an elliptic analogue of the result of
\cite{DF}.

\section{About the proof}

\subsection{A rational analogue}

In \cite{Rat}, we constructed a Hopf algebra cocycle in the Yangian
double $A = DY(\SL_{2})$,
conjugating Drinfeld's coproduct to the usual one. $DY(\SL_{2})$ is a
rational version of $U_{\hbar}\G(\tau)$ (it was first introduced in
\cite{Kh}). This means that it has the same
presentation as $U_{\hbar}\G(\tau)$, replacing the theta-functions by
their arguments. 

For $x = e,f,h$, let $x_{n}$ be
the analogue of $x[z^{n}]$. $DY(\SL_{2})$ contains two
``negative modes'' and ``nonnegative modes'' subalgebras $A^{<0}$ and
$A^{\ge 0}$, repectively
generated by $K$ and the $x_{n},n<0$ and $D$ and the 
$x_{n},n\ge 0$, $x = e,f,h$. 

The Yangian coproduct $\De_{\Yg}$ on 
$DY(\SL_{2})$ defines a Hopf algebra structure
on $DY(\SL_{2})$, for which $A^{\ge 0}$ and
$A^{<0}$ are both Hopf subalgebras. On the other hand, the rational
analogues $\De_{\rat}$ and $\bar\De_{\rat}$ of $\De$ and $\bar\De$
have the following properties: 
$$
\Delta_{\rat}(A^{\geq 0}) \subset A \otimes A^{\geq 0}, \quad
\Delta_{\rat}(A^{<0}) \subset A^{<0} \otimes A, 
$$
$$
\bar \Delta_{\rat}(A^{\geq 0}) \subset A^{\geq 0} \otimes A, \quad
\bar \Delta_{\rat}(A^{<0}) \subset A \otimes A^{<0}.  
$$

Based on the study of Hopf algebra duality within $DY(\SL_{2})$, we show
that we have a decomposition of the rational analogue of $F$, 
$F_{\rat}$, as a product  $F_{\rat} = F_{2} F_{1}$, with 
$F_{1}\in A^{<0}
\otimes A^{\ge 0}$ and 
$F_{2} \in A^{\ge 0}
\otimes A^{< 0}$. 

Then the twist $\De = F_{1}\De_{\rat}F_{1}^{-1}$ coincides with the
conjugation $ F_{2}^{-1}\bar\De_{\rat}F_{2}$. It follows that $\De$
satisfies both 
$$
\De(A^{\ge 0}) \subset A^{\ge 0}\otimes A^{\ge 0} \quad 
\on{and} \quad 
\De(A^{<0}) \subset A^{<0}\otimes A^{<0}. 
$$

On the other hand, $\De$ defines a Hopf algebra structure on
$DY(\SL_{2})$. Indeed, the associator corresponding to $F_{1}$ is
expressed as 
$$
\Phi = F_{1}^{(12)}(\Delta\otimes 1)(F_{1})
\left( F_{1}^{(23)}(1\otimes \Delta)(F_{1})\right)^{-1} .  
$$
We have clearly $\Phi \in A^{<0} \otimes A \otimes A^{\geq 0}$. 
Since we also have 
$$
\Phi = \left( (\bar\Delta\otimes 1)(F_{2}) F_{2}^{(12)} \right)^{-1}
(1\otimes \bar\Delta)(F_{2}) F_{2}^{(23)} ,
$$
we also see that $\Phi\in A^{\geq 0} \otimes A \otimes
A^{<0}$. Therefore, 
$\Phi = 1 \otimes a \otimes 1$, for a certain $a\in A$. On the 
other hand, as $\Phi$ is obtained by twisting a quasi-Hopf structure, it 
should satisfy the compatibility condition (see \cite{D-qH})
$$
(\Delta \otimes id \otimes id)(\Phi)
(id \otimes id \otimes \Delta)(\Phi)
=
(\Phi\otimes 1)(id \otimes \Delta\otimes id)(\Phi) 
(1\otimes \Phi) ,
$$
so that $a=1$. 
Therefore, $\Phi = 1$ and $F_{1}$ satisfies the Hopf cocycle
equation. It follows that $\De$ defines Hopf algebra structure on
$DY(\SL_{2})$, which we can show to coincide with $\De_{\Yg}$ (see
\cite{Rat}).  

The strategy followed in \cite{EF} can be described as
follows. $U_{\hbar}\G(\tau)$ contains an analogue of $A^{\ge 0}$, its
subalgebra $U_{\hbar}\G_{\cO}$ generated by the $x[\eps],\eps\in\cO,x =
e,f,h$. It contains no analogue of $A^{<0}$. However, there are certain
subalgebras of  ${\on{Hol}}(\CC - \Gamma , U_{\hbar}\G(\tau)^{\otimes 2})$  
and
${\on{Hol}}(\CC - \Gamma, U_{\hbar}\G(\tau)^{\otimes 3})$, playing the
roles of $A^{\ge 0} \otimes A^{<0}$ and $A^{<0} \otimes A^{\ge 0}$,
respectively $A^{\ge 0, <0} \otimes A^{\otimes 2}$ and $A^{\otimes 2}
\otimes A^{\ge 0,<0}$. We will decompose $F$ as a product of elements of
these algebras. This will give rise to a solution of the twisted cocycle
equation, and by BBB, to a solution of the DYBE in ${\on{Hol}}(\CC -
\Gamma , U_{\hbar}\G(\tau)^{\otimes 2})$. This solution will happen to
coincide with $R(z,\la)$ in suitable representations. Therefore, the
corresponding $L$-operators will satisfy the elliptic quantum group
relations. Let us see now in more detail how this program is realized.  

\subsection{Subalgebras of 
$Hol(\CC - \Gamma, U_{\hbar}\G (\tau)^{\otimes 2}) $  
and
$Hol(\CC - \Gamma, U_{\hbar}\G (\tau)^{\otimes 3})$}

For $X$ a vector space, call $\Hol(\CC - \Gamma, X)$ the space 
of holomorphic functions from $\CC - \Gamma$ to $X$ and set $\Hol(\CC-
\Gamma) = 
\Hol(\CC - \Gamma,\CC)$. The parameter in $\CC - \Gamma$ will be
identified with the spectral parameter $\la$. 

\begin{defin}
Let us define $A^{-+}$ to be 
the subalgebra of $\Hol(\CC - \Gamma , A(\tau)^{\otimes 2})$
generated (over $\Hol(\CC - \Gamma)$) by $h^{(2)}$ and the 
$e^{-(1)}_{-\la -\hbar h^{(2)}}[\eps] f^{(2)}[r]$, with $\eps\in \cK$
and $r\in \cO$, and $A^{+-}$ as the subalgebra of $\Hol(\CC - \Gamma,
A(\tau)^{\otimes 2})$ 
generated (over $\Hol(\CC - \Gamma)$) 
by $h^{(2)}$ and the $e^{(1)}[r]f^{-(2)}_{\la +\hbar h^{(2)} 
-2\hbar}[\eps]$, with $r\in\cO$, $\eps\in \cK$. 
\end{defin}

\begin{defin}
Let us define $A^{-,\cdot,\cdot}$ as the subspace of the algebra
$\Hol(\CC - \Gamma, A(\tau)^{\otimes 3})$, linearly 
spanned (over $\Hol(\CC - \Gamma)$) by the elements of the form 
\begin{equation} \label{medved}
\xi' = e^{-(1)}_{-\la -\hbar ( h^{(2)} + h^{(3)}) 
}[\eta_{1}]\cdots
e^{-(1)}_{-\la -\hbar ( h^{(2)} + h^{(3)}) -2(n-1)\hbar
}[\eta_{n}] 
(1\otimes a \otimes b), 
\end{equation}
$n\ge 0$ (recall that the empty product is equal to $1$), 
where $\eta_{i}\in \cK$, and $a,b\in A(\tau)$ are 
such that $[h^{(1)} + h^{(2)} + h^{(3)}, \xi'] =0$; 
and $A^{\cdot,\cdot,+}$ as the subspace of $\Hol(\CC - \Gamma, 
A(\tau)^{\otimes 3})$ 
spanned (over $\Hol(\CC - \Gamma)$) by the elements of the form 
$$
\eta' = ( a^{\prime}\otimes b' \otimes 1) 
f^{(3)}[r_{1}]\cdots
f^{(3)}[r_{n}] (h^{(3)})^{s}, \quad n,s\ge 0, 
$$
where $a',b'\in A(\tau)$, $r_{i}\in \cO$, 
and such that $[h^{(1)} + h^{(2)} + h^{(3)}, \eta'] =0$.
\end{defin}

\begin{defin}
$A^{+,\cdot,\cdot}$ is the subspace of the algebra $\Hol(\CC - \Gamma,
A(\tau)^{\otimes 3})$ linearly spanned (over $\Hol(\CC - \Gamma)$)
by the elements of the form
\begin{equation} \label{xi'}
\xi' = e^{(1)}[r_{1}] \ldots e^{(1)}[r_{n}] (1\otimes a \otimes b), 
\quad n\ge 0, 
\end{equation}
where $r_{i}\in \cO$, and $a,b\in A(\tau)$ are such that $[h^{(1)} + 
h^{(2)} + h^{(3)}, \xi']=0$.  
 
$A^{\cdot,\cdot,-}$ is the subspace of $\Hol(\CC - \Gamma, 
A(\tau)^{\otimes 3})$ 
linearly spanned (over $\Hol(\CC - \Gamma)$) by the elements of the form
\begin{equation} \label{eta'}
\eta' = (a'\otimes b'\otimes 1)
f^{-(3)}_{\la +\hbar h^{(3)} -2\hbar}[\eta_{1}] \ldots 
f^{-(3)}_{\la +\hbar h^{(3)} -2\hbar}[\eta_{n}] (h^{(3)})^{s} , 
\quad n,s\ge 0, 
\end{equation}
where $\eta_{i}\in \cK$, and $a',b'\in A(\tau)$ are such that $[h^{(1)} + 
h^{(2)} + h^{(3)}, \eta']=0$.  
\end{defin}

\begin{prop} \label{egill}
$A^{-,\cdot,\cdot}$, $A^{\cdot,\cdot,+}$, $A^{+,\cdot,\cdot}$ and
$A^{\cdot,\cdot,-}$ are subalgebras
of $\Hol(\CC - \Gamma, A(\tau)^{\otimes 3})$. We have
\begin{equation} \label{genet}
(\Delta\otimes 1)(A^{-+}) \subset A^{-,\cdot,\cdot}
\cap A^{\cdot,\cdot,+} , 
\quad
(1\otimes \Delta)(A^{-+}) \subset A^{-,\cdot,\cdot} 
\cap A^{\cdot,\cdot,+},   
\end{equation}
\begin{equation} \label{tawil}
(\bar\Delta \otimes 1)(A^{+-}) \subset A^{+,\cdot,\cdot}
\cap A^{\cdot,\cdot,-}, 
\quad
(1\otimes \bar\Delta)(A^{+-}) \subset A^{+,\cdot,\cdot,} \cap 
A^{\cdot,\cdot,-}.  
\end{equation}
\end{prop}

The first statement is a consequence of the following relations between
half-currents:

\begin{lemma}
The generating series $e^{\pm}_{\la}(z), f^{\pm}_{\la }(z)$  
satisfy the following relations: 
\begin{align} \label{eeps:eeps'}
{ {\theta(z-w-\hbar)} \over {\theta(z-w)} } & e^{\eps}_{\la+\hbar}(z)  
e^{\eps'}_{\la-\hbar}(w)
+
\eps\eps'
{{\theta(w-z-\la)\theta(-\hbar)} \over {\theta(w-z)\theta(-\la)} } 
e^{\eps'}_{\la+\hbar}(w) e^{\eps'}_{\la-\hbar}(w)
\\ \nonumber & =
{ {\theta(z-w+\hbar)} \over {\theta(z-w)} } e^{\eps'}_{\la+\hbar}(w) 
e^{\eps}_{\la-\hbar}(z)
+
\eps\eps'
{ {\theta(z-w-\la)\theta(-\hbar)} \over {\theta(z-w)\theta(-\la)} } 
e^{\eps}_{\la+\hbar}(z) e^{\eps}_{\la-\hbar}(z),
\end{align}
\begin{align} \label{feps:feps'}
{ {\theta(z-w+\hbar)} \over {\theta(z-w)} } & f^{\eps}_{\la-\hbar}(z)  
f^{\eps'}_{\la+\hbar}(w)
+
\eps\eps'{ {\theta(w-z-\la)\theta(\hbar)} \over {\theta(w-z)\theta(-\la)} } 
f^{\eps'}_{\la-\hbar}(w) f^{\eps'}_{\la+\hbar}(w)
\\ \nonumber & =
{ {\theta(z-w-\hbar)} \over {\theta(z-w)} } f^{\eps'}_{\la-\hbar}(w) 
f^{\eps}_{\la+\hbar}(z)
+
\eps\eps'
{ {\theta(z-w-\la)\theta(\hbar)} \over {\theta(z-w)\theta(-\la)} } 
f^{\eps}_{\la-\hbar}(z) f^{\eps}_{\la+\hbar}(z) ,  
\end{align}
where $\eps,\eps'$ take the values $+,-$. 
\end{lemma}

To prove the statement about coproducts in Prop. \ref{egill}, it is then
sufficient to prove
it from generators of $A^{+-}$ and $A^{-+}$ (see \cite{EF}, Lemmas
1.2-4). 

A useful argument in sect. 4.1 was that the intersection of 
$A^{\ge 0} \otimes A \otimes A^{<0}$ and $A^{<0} \otimes A \otimes
A^{\ge 0}$ is reduced to $1 \otimes A \otimes 1$. Similarly, we have: 

\begin{prop} \label{inters}
We have 
$$
A^{+,\cdot,\cdot}\cap A^{-,\cdot,\cdot} = \Hol(\CC - \Gamma , 1\otimes 
(A(\tau)^{\otimes2})^{\HH}), 
\ 
A^{\cdot,\cdot,+}\cap A^{\cdot,\cdot,-} = \Hol(\CC - \Gamma ,
(A(\tau)^{\otimes 2})^{\HH}
\otimes \CC[h]),  
$$
where $(A(\tau)^{\otimes 2})^{\HH}$ are the elements of 
$A(\tau)^{\otimes 2}$ commuting with $h^{(1)} + h^{(2)}$. 
\end{prop}

\subsection{Decomposition of $F$}

We then have:

\begin{prop} \label{dyn:decomp}
There is a unique a decomposition of $F$ as 
\begin{equation} \label{decomp:ell}
F= F^{2}_{\la}F^{1}_{\la} , \quad \on{with} \quad  
F^{1}_{\la}\in A^{-+}\quad \on{and} \quad 
F^{2}_{\la}\in A^{+-},  
\end{equation}
with $(\varepsilon\otimes 1)(F_{i}) = (1\otimes \varepsilon)(F_{i}) = 1$, 
$i=1,2$. 
\end{prop}

As in the rational case, the proof of this fact relies on a study of the
Hopf duality between the opposite quantum nilpotent current subalgebras
of $U_{\hbar}\G(\tau)$. 

\begin{lemma} \label{approx}
We have an expansion 
$$
F^{1}_{\la} \in 1 + \hbar \sum_{i\ge 0} e^{-(1)}_{-\la -\hbar h^{(2)}}
[e_{i;-\la}] f^{(2)}[e^{i}]+ U_{\hbar}\N_{+}(\tau)^{\ge 2}
\otimes U_{\hbar}\N_{-}(\tau)^{\ge 2}\CC[h], 
$$
where $U_{\hbar}\N_{\pm}(\tau)^{\ge 2} = \oplus_{i\ge 2}
U_{\hbar}\N_{\pm}(\tau)^{[i]}$, and  $U_{\hbar}\N_{\pm}(\tau)^{[i]}$ is
the linear span of the products of $i$ terms of the form $e[\eps]$ in
$\pm = +$, and $f[\eps]$ if $\pm = -$. 
\end{lemma}

\subsection{Twisted cocycle equation}

\begin{prop} \label{tw-cocycle}
The family $(F^{1}_{\la})_{\la\in \CC - \Gamma}$ satisfies the twisted
cocycle equation 
\begin{equation} 
F^{1(12)}_{\la +\hbar h^{(3)}} (\Delta\otimes 1)(F^{1}_{\la})
=
F^{1(23)}_{\la} (1\otimes \Delta)(F^{1}_{\la}).
\end{equation}  
\end{prop}

To show this result, we proceed as follows. Define for $\la\in \CC -
\Gamma$,  
$$
\Phi_{\la} = F_{\la+\hbar h^{(3)}}^{1(12)}
(\Delta\otimes 1)(F_{\la}^{1})
\left( F_{\la}^{1(23)} (1\otimes \Delta)(F_{\la}^{1})\right)^{-1} . 
$$
Then (\ref{genet}) implies that $(\Phi_{\la})_{\la\in\CC - \Gamma}$ 
belongs to $ A^{-,\cdot,\cdot} \cap A^{\cdot,\cdot,+}$. On the other
hand, using (\ref{F:cocycle}), we may rewrite $\Phi_{\la}$ as 
$$
\Phi_{\la} = \left((\bar\Delta\otimes 1)(F^{2}_{\la})(F^{2(12)}_{\la +
\hbar h^{(3)}}) \right)^{-1} (1\otimes
\bar\Delta)(F_{\la}^{2}) 
(1\otimes F^{2(23)}_{\la}). 
$$
Then (\ref{tawil}) then implies that $(\Phi_{\la})_{\la\in\CC - \Gamma}$
belongs to $A^{+,\cdot,\cdot} \cap A^{\cdot,\cdot,-}$. 

Prop. \ref{inters} then implies that 

\begin{equation} \label{form-of-Phi}
\Phi_{\la}  = \sum_{i\ge 0}1 \otimes a^{(i)}_{\la} \otimes h^{i}. 
\end{equation}
Define now for $x\in A(\tau)$, 
$$
\Delta_{\la}(x) = F^{1}_{\la} \Delta(x) (F^{1}_{\la})^{-1} ; 
$$
$(\Delta_{\la})_{\la\in\CC - \Gamma}$ and $(\Phi_{\la})_{\la\in\CC - \Gamma}$ 
satisfy the compatibility condition
\begin{align} \label{maoz}
& (\Delta_{\la +\hbar (h^{(3)} + h^{(4)})} \otimes 1 \otimes 1)
(\Phi_{\la})
(1\otimes 1 \otimes \Delta_{\la}) (\Phi_{\la})
\\ & \nonumber
=
\Phi^{(123)}_{\la +\hbar h^{(4)}} (1\otimes \Delta_{\la +\hbar h^{(4)}}
\otimes 1)
(\Phi_{\la}) \Phi_{\la}^{(234)}. 
\end{align}

It follows that $\Phi_{\la} = 1$.

\subsection{The BBB result}

The paper \cite{BBB}, sect. 2, contains the following result:  

\begin{prop} (see \cite{BBB}) \label{bbb}
Let $(\cA,\Delta_{\infty}^{\cA},
\cR_{\infty}^{\cA})$ 
be a quasi-triangular Hopf algebra, with a fixed element $\wt h$. 
Let $F(\la)$ be a family of invertible elements of $\cA \otimes \cA$, 
parametrized by some subset $U \subset \CC$. 
Set $\Delta (\la) = \Ad(F(\la)) \circ 
\Delta^{\cA}_{\infty}$.  Suppose that the identity
\begin{equation} \label{tw-coc}
F^{(12)}(\la+\hbar \wt h^{(3)}) (\Delta_{\infty}^{\cA}\otimes
1)(F(\la))
=
F^{(23)}(\la) (1\otimes \Delta_{\infty}^{\cA})(F(\la))
\end{equation}
is satisfied. Then we have 
\begin{equation} \label{tw:coass}
(\Delta (\la +\hbar \wt h^{(3)}) \otimes 1) \circ \Delta (\la)
=
(1 \otimes \Delta (\la)) \circ \Delta (\la) ,
\end{equation}
and  if we set $\cR(\la) = F^{(21)}(\la)\cR_{\infty}^{\cA}F(\la)^{-1}$, 
we have the identity
\begin{equation} \label{DYBE}
\cR^{(12)}(\la) \cR^{(13)}(\la +\hbar \wt h^{(2)}) \cR^{(23)}(\la) 
=
\cR^{(23)}(\la +\hbar \wt h^{(1)}) \cR^{(13)}(\la) \cR^{(12)}(\la +
\hbar \wt h^{(3)}) . 
\end{equation}
\end{prop}

Equation (\ref{DYBE}) is the algebraic version of the DYBE. 

\subsection{End of the proof}

Apply now Prop. \ref{bbb} to $\cA = U_{\hbar}\G(\tau)$, 
$\De_{\infty}^{\cA} = \De$, $\cR_{\infty}^{\cA} = \cR_{\infty}$ and
$F(\la) = F_{\la}^{1}$. Set  $\cR_{\la} = (F^{1}_{\la})^{(21)} \cR_{\infty}
(F^{1}_{\la})^{-1}$. We have: 

\begin{corollary} The family
$(\cR_{\la})_{\la\in \CC - \Gamma}$ satisfies the dynamical Yang-Baxter
relation 
\begin{equation} \label{DYBE:ell}
\cR^{(12)}_{\la} \cR^{(13)}_{\la +\hbar h^{(2)}} \cR^{(23)}_{\la} 
=
\cR^{(23)}_{\la +\hbar h^{(1)}} \cR^{(13)}_{\la} \cR^{(12)}_{\la +
\hbar h^{(3)}}. 
\end{equation}
\end{corollary}

The procedure to derive $RLL$ relations from a Yang-Baxter-like equation
is to first study finite-dimensional representations of the algebra, and
then take the image of the YBE in these representations (see
\cite{RS}).

\begin{prop}  \label{fd:repres} (see \cite{HGQG}, Prop. 9)
There is a morphism of algebras 
$\pi_{\zeta} : A(\tau) \to \End(\CC^{2}) \otimes \cK_{\zeta}[\pa_{\zeta}]
[[\hbar]]$, defined by the formulas 
$$
\pi_{\zeta}(K) = 0, \quad \pi_{\zeta}(D) = \Id_{\CC^{2}} \otimes \pa_{\zeta}, 
$$
$$
\pi_{\zeta}(h[r]) = 
E_{11}  \otimes \left( {{2}\over{1+q^{\pa}}} r \right)
(\zeta) 
- E_{-1-1} \otimes \left( {{2} \over {1+q^{-\pa}}} r \right) (\zeta) ,
\quad
r\in \cO, 
$$
$$
\pi_{\zeta}(h[\la]) = E_{11} \otimes \left( {{ 1 - q^{-\pa}}\over{\hbar
\pa }} \la \right) (\zeta) - E_{-1-1} \otimes 
\left( {{ q^{\pa} - 1} \over {\hbar\pa}} \la \right) (\zeta) , 
\quad \la\in L_{0}, 
$$
$$
\pi_{\zeta}(e[\eps]) = { { \theta(\hbar)} \over \hbar } 
E_{1,-1}\otimes \eps(\zeta), 
\quad 
\pi_{\zeta}(f[\eps]) = E_{-1,1}\otimes \eps(\zeta), 
\quad \eps\in \cK.  
$$
\end{prop}

The image of $\cR_{\la}$ by these representations is computed as
follows: 

\begin{lemma} \label{adar}
The image of $\cR_{\la}$ by $\pi_{\zeta} \otimes \pi_{\zeta'}$ is 
\begin{equation} \label{im-R}
(\pi_{\zeta} \otimes \pi_{\zeta'})(\cR_{\la+\hbar}) = 
A(\zeta,\zeta') R^{-}( \zeta-\zeta', \la) , 
\end{equation}
where $R^{-}(z,\la)$ has been defined in (\ref{R-mat}). 
\end{lemma}

The computation relies on Lemma \ref{approx}. 

Appling $id \otimes \pi_{\zeta} \otimes
\pi_{\zeta'}$, $\pi_{\zeta} \otimes \pi_{\zeta'} \otimes id$ and
$\pi_{\zeta} \otimes id \otimes \pi_{\zeta'}$ to (\ref{DYBE:ell}), after
the change of $\la$ into $\la+\hbar$, we find (\ref{Lpm:Lpm}) and
(\ref{L+:L-}).

\end{document}